\documentclass[journal]{IEEEtran}
\usepackage{graphicx} % Required for inserting images
\usepackage{stackengine}

\usepackage{amssymb}
\usepackage{svg}
\usepackage{algorithm}
\usepackage[noend]{algpseudocode}
\usepackage{multirow}
\usepackage{subfigure}
\usepackage{multirow}
\usepackage{array}
\usepackage{booktabs}
\usepackage{subcaption}
\usepackage[font=small]{caption}
\usepackage{amsmath}
\usepackage{stmaryrd}
\usepackage{mathtools}
\DeclareMathOperator*{\argmin}{argmin}
\algblockdefx[Ifm]{Ifm}{EndIfm}[0]{\textbf{if} $( $}{$)$}
\algcblockdefx[Then]{Ifm}{Then}{EndThen}{$\textbf{then}$}{}
\algcblockdefx[Elsem]{Then}{Elsem}{EndElsem}{$\}$ \textbf{else} $\{$}{$\}$}

\newcommand\barbelow[1]{\stackunder[1.2pt]{$#1$}{\rule{.8ex}{.075ex}}}

\ifCLASSINFOpdf
\else
   \usepackage[dvips]{graphicx}
\fi
\usepackage{url}

\usepackage{graphicx}

\begin{document}
\title{Mathematical Optimization-Based Period Estimation with Outliers and Missing Observations}

\author{Romain Puech, Vincent Gouldieff %\IEEEmembership{Member, IEEE}
\thanks{This work was in part supported by Safran Data Systems, Les Ulis, France. August 2023.
Pre-publication version. Copyright statement: Copyright may be transferred without notice, after which this version may no longer be accessible.}
\thanks{Romain Puech is associated with École polytechnique, Palaiseau, France (email: romain.puech@alumni.polytechnique.org).}
\thanks{ Vincent Gouldieff is associated with Safran Data Systems, Les Ulis, France.}}

\markboth{August 2023}
{Shell \MakeLowercase{\textit{et al.}}: Bare Demo of IEEEtran.cls for IEEE Journals}
\maketitle

\begin{abstract}
We consider the frequency estimation of periodic signals using noisy time-of-arrival (TOA) information with missing (sparse) data contaminated with outliers. We tackle the problem from a mathematical optimization standpoint, formulating it as a linear regression with an unknown increasing integer independent variable and outliers. Assuming an upper bound on the variance of the noise, we derive an online, parallelizable, near-CRLB optimization-based algorithm amortized to a linear complexity. We demonstrate the outstanding robustness of our algorithm to noise and outliers by testing it against diverse randomly generated signals. Our algorithm handles outliers by design and yields precise estimations even with up to 20\% of contaminated data.
\end{abstract}


\IEEEpeerreviewmaketitle

\section{Introduction}
\IEEEPARstart{T}{he} precise identification of the frequency represents a cornerstone task of periodic signal processing. Pulse periodic trains constitute a subset of signals that have emerged as a subject of keen interest due to their applicability in diverse domains such as radar systems \cite{radar_application}\cite{radar_recent}, wireless communications \cite{comunication}, and medical instrumentation \cite{ECG}\cite{neural}.
Moreover, they can be used to describe more complex periodic signals \cite{periodic_pulse_trains_are_everywhere} given appropriate preprocessing, using zero-crossing \cite{zero_crossing}\cite{zero_crossing_modern} or a correlator for example.
For such trains, the time-of-arrival (TOA) $y_i$ of the $i^{th}$ pulse is related to the period $T$, phase $b$ and error (noise) $\epsilon_i$ by the linear relationship:
\begin{equation}\label{equation:linear}
    y_i = x_i \times T + b + \epsilon_i,\ \text{with } x_i = i
\end{equation}

Nevertheless, the estimation of pulse repetition interval or period in pulse trains becomes challenging in the eventuality of missing pulses ($x_i$ becomes unknown), and when the data is contaminated with false pulse detections (outliers).

\begin{figure}[htbp]\label{fig:subfig1}
    \centering
        \includegraphics[width=0.489\textwidth]
        {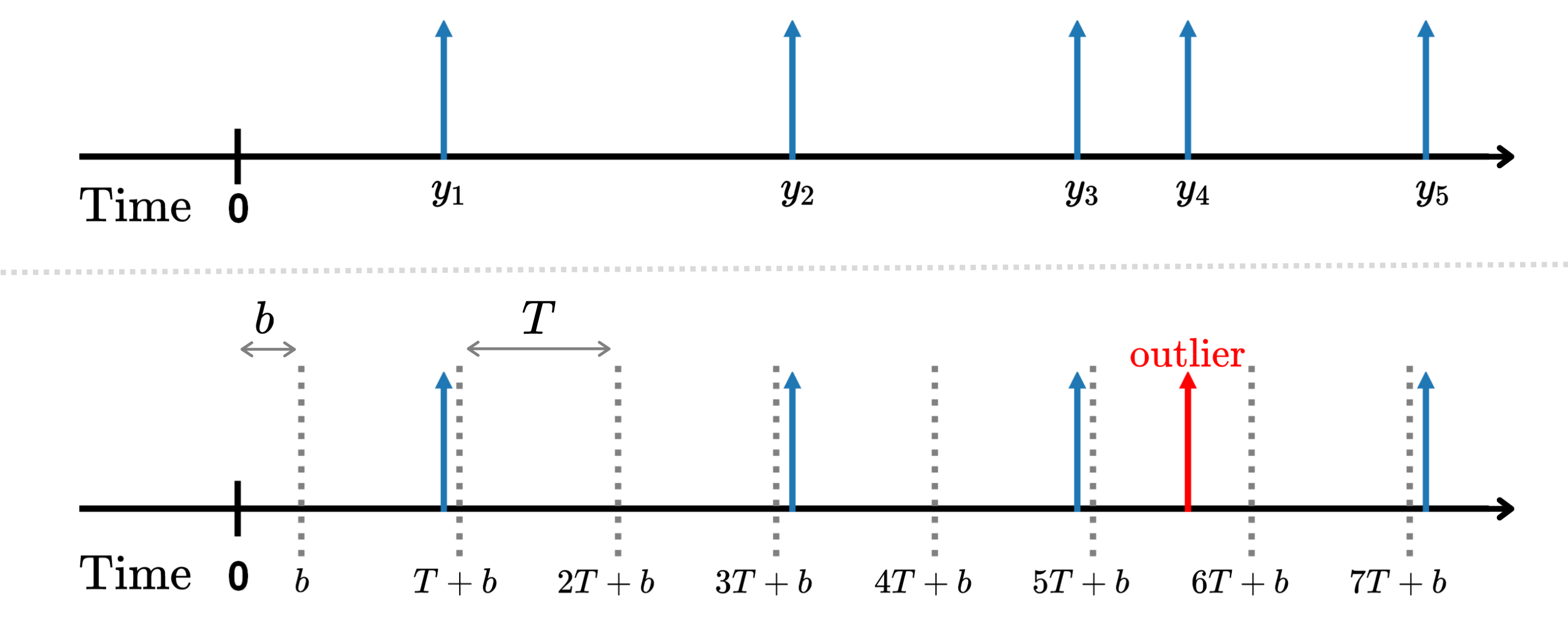}
        %\label{fig:repr}
    \caption{Above: the information given as input. Below: the pulse train that generated the input, and that we want to find back}
    \label{fig:repr}
\end{figure}
%Such noise level can for example happen in ambulatory electrocardiograms \cite{ECG} leading to poor performances of beat detection algorithms.

%Ranney and Tom \cite{survey_of_methods} present various algorithms based on spectral analyses notably via auto-correlation, or using histogram methods. However, these methods are not tailored for missing pulses and misdetections. 
Sadler and Casey \cite{periodic_pulse_signal_analysis} address this problem, using a modified Euclidean algorithm taking as input the TOAs. However, in the presence of outliers, their algorithm is only effective if most of the pulses are non-missing, that is, if the duration between two consecutive pulses is typically of one period.
Later, using a geometric interpretation of the problem, Sidiropoulos \textit{et al.} \cite{quasi_ml} derived a quasi-maximum likelihood approach called separable least squares line search (SLS2) achieving the Cramér-Rao lower bound (CRLB). A review of different algorithms is proposed in \cite{improved_quasi_ml}, which features an improved SLS2 algorithm and performs extensive empirical comparisons. None of the algorithms presented addresses by design the question of outliers, and most display sub-optimal precision when the ratio between the standard deviation of the noise and the period exceeds $1/50$. Moreover, the tests are performed when few pulses are missing (generated by a Bernoulli process of parameter $p=0.5$). %An ML estimator assuming Gaussian noise and a lower bound on the period is derived in the same paper, although impractical due to its exponential complexity. 
A fast log-linear estimator is derived in \cite{fast_but_not_precise}, but at the expense of some accuracy by quantizing the observations. A linear estimator of comparable precision is presented in \cite{linear}.

%% Figure was there before

Our paper aims to address the intricate issue of period estimation in the eventuality of numerous missing pulses and taking into account outliers. A visual representation of our problem is displayed in Fig.\ref{fig:repr}. We adopt a mathematical optimization perspective on the problem by formulating and solving it as a non-convex mixed-integer non-linear optimization problem, taking as input the time-of-arrival (TOA) information and an upper bound $\overline{\sigma}$ on the variance of its noise. We express the period estimation as a linear regression with an unknown increasing integer independent variable and outliers and use a branch-and-bound-inspired algorithm to solve it. 

Our contributions are summarized as follows:
\begin{itemize}
    \item We formulate the period estimation problem as an optimization problem and exhibit a tree representation for its search space.
    \item We present an online, parallelizable, near-CRLB estimation algorithm for the period, displaying outstanding resilience to outliers and noise, amortized to a linear complexity in the number of pulses.
    %\item We provide a proof of the algorithm.
    \item We test our algorithm on generated data.
\end{itemize}
In section \ref{section:modelisation} we detail our representation of the problem as a mathematical optimization one and present our algorithm and its proof in section \ref{section:resolution}. We test our algorithm on generated signals in section \ref{section:test}.

\section{Modelisation}\label{section:modelisation}
We define the problem as a mathematical optimization one.
\paragraph{Input}
Consider a periodic signal transformed into a single periodic pulse train with period $T$. The TOA for a fraction of these pulses is extracted up to an error $\epsilon \sim \mathcal{N}(0,\sigma^2)$, using a matched filter or a similar preprocessing algorithm. Since $\sigma$ mainly depends on the preprocessing accuracy, we assume a known upper bound $\overline{\sigma}$ for it. We set $\overline{\sigma}/T < 1/10$ as a reasonable assumption so that the order of the pulses is preserved with probability close to $1$.
The input $\boldsymbol{y}= (y_i)_{1 \leq i \leq n}$ consists of $n$ detected TOAs, including outliers.

\paragraph{Unknowns}
We introduce $\boldsymbol{x} \in (\mathbb{N}\cup\{-1\})^n$, where $x_i$ represents the period index of the $i^{th}$ detected pulse (with $-1$ for a false detection). For example, $x_3 = 5$ means that the third detected pulse started in period 5. For the sake of clarity, we assume that the first two points are not outliers and set $x_1 = 0$ by subtracting $y_1$ from each $y_i$. We will easily drop that assumption in section \ref{section:resolution}. We denote $b \in [0,T)$ the phase and $I$ the set of indices of non-outlier pulses. This way, $i \in I$ implies that $y_i$ is the TOA of a real pulse, as opposed to an outlier. Similarly to (\ref{equation:linear}) we have:
\begin{equation}
    \forall{i \in I, \ }y_i = x_i \times T + b + \epsilon_i
\end{equation}

\paragraph{Optimization problem}
If $\boldsymbol{x}$ and $I$ were known, finding $T$ and $b$ would amount to a linear regression, achieving the CRLB with an ordinary least squares (OLS) estimator as stated by the Gauss-Markov theorem. The intricacy of the problem comes from the fact that
none of $\boldsymbol{x}$, $T$, $b$ or $I$ is known. Excluding outliers, the problem can thus be formulated as a simple linear regression with an unknown increasing integer independent variable.

We introduce $N \in \mathbb{N}$ as an upper bound on the number of consecutive undetected pulses, and $\mathcal{X}$ as the set of valid sequences for $\boldsymbol{x}$. We call $\boldsymbol{x}$ valid if it is an $n$-tuple of $m\leq n$ increasing non-negative integers, with a difference of at most N between two consecutive items, to which we insert $n-m$ times $\{-1\}$ at any indices to represent the outliers.
Finding the period $T$ then amounts to solving the following least squares problem:

\begin{equation}\label{objective}
    \argmin\limits_{\substack{x \in \mathcal{X} \\ T \in \mathbb{R}_{>0} \\ b \in [0,T)}} \sum\limits_{i \in I}(y_i - (T \times x_i + b))^2
\end{equation}
The problem is easily seen to be non-convex: if $(x,T,b)$ is a solution, so is $(2x,T/2,b)$. Another difficulty to take into account is that the minimum is taken over all $i\in I$ and not $i \in \{1,\ldots,n+m\}$: the set of variables we are optimizing on is also unknown, implying that solving the optimization over all the input won't necessarily lead to the sought solution: we also need a way to flag outliers. The formulation can be simplified by noticing that for any $\boldsymbol{x}$ given, $T$ and $b$ can only be the OLS estimator associated with $(x_i)_{i\in I}$ in order to reach optimality. Indeed, the objective function is convex in $T$, so it admits a global minimizer, that coincides with the OLS.   

Putting aside outliers for a time, one can consider the naive approach of trying all possible $\boldsymbol{x}$ in the search space $\mathcal{X}$ and outputting one leading to the lowest value of the objective function in (\ref{objective}). The cost of computing the objective function being linear in $n$, the complexity of the brute force algorithm is 
$O(n|\mathcal{X}|)$. Recalling the upper bound $N$ on the distance between adjacent elements of $\boldsymbol{x}$, we get an exponential complexity of $O(nN^{n-1})$.
%%%%%%%%%%%%%%%%%%%%%%%%%%%%%

 \section{resolution}\label{section:resolution}
 While a brute force approach would be intractable, we will focus on representing the set of solutions with a tree structure and traverse it in a branch-and-bound manner to avoid branches leading to a suboptimal solution at a specified confidence level $c$, and flagging spurious pulse detections. This approach guarantees an optimal solution at a chosen confidence level $c$ representing a trade-off between speed and robustness. %Anticipating on the testing section, this algorithm allowed us to solve problems with $n=1000$ and $N=100$ in a fraction of a second. 
 The complexity of our approach is amortized to linear in $n$ for $N$ fixed, massively improving on the exponential complexity of the naive approach.
 %FLOW CHART WAS THERE%
 \paragraph{Tree representation}  As previously introduced, we denote by $\mathcal{X}$ the search space. For $1 \leq k \leq n$, we denote by $\mathcal{X}_{(x_1, \ldots, x_k)} \subset \mathcal{X}$ the set of all possible solutions $\boldsymbol{x}$ beginning with $(x_1, \ldots, x_k)$. This way, $\mathcal{X}$ can be represented by a tree with $\mathcal{X}_{(0)}$ as a root, $\mathcal{X}_{(0,1)},\ldots,\mathcal{X}_{(0,N)},\mathcal{X}_{(0,-1)}$ at depth one, and so on, with solutions specified up to the $k^{th}$ element at depth $k$. Each leaf of this tree represents a unique potential solution. Our algorithm consists of exploring this tree structure using a depth-first traversal while pruning the appropriate branches as early as possible. An accurate evaluation of how unlikely is a certain branch to lead to an optimal solution is thus crucial to the efficiency of our algorithm.
 %and is discussed in paragraph \ref{tests}.
 \begin{figure}
     \centering
     \includegraphics[width=0.62\columnwidth]{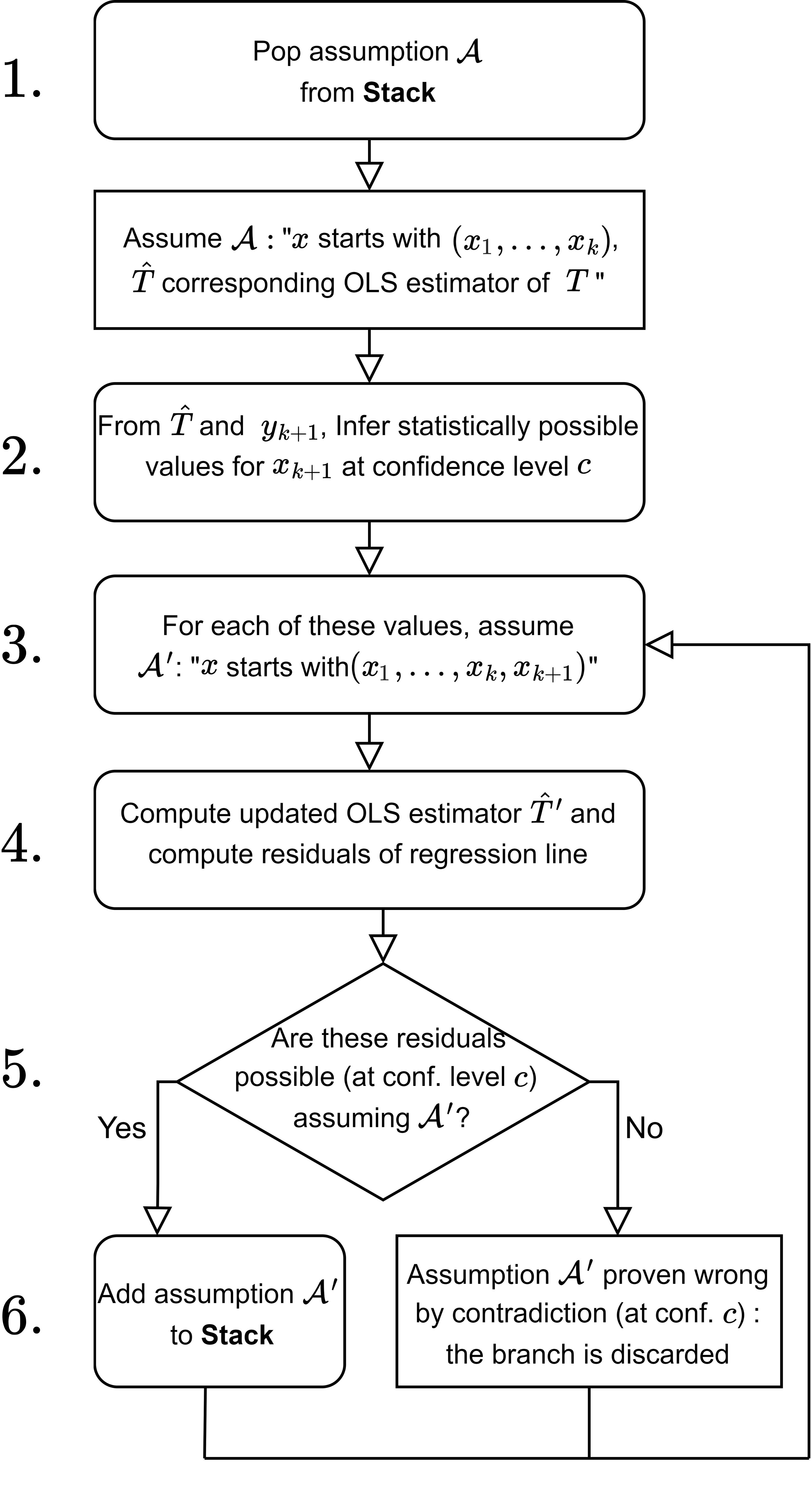}
     \caption{Single iteration of the algorithm.}\label{fig:flow}
 \end{figure}
\paragraph{Iteration on a branch}
For a branch $\mathcal{X}_{(x_1,\ldots,x_k)}$, the working principle of the algorithm is displayed in Fig. \ref{fig:flow} and in Algorithm \ref{alg:iteration}. As represented, the sub-branches stemming from $\mathcal{X}_{(x_1,\ldots,x_k)}$ are considered one by one. These are first obtained in step 2 in Fig. \ref{fig:flow} (detailed in the next paragraph) using the estimate $\hat{T}$ of the parent branch. For each of them, we then fit an updated OLS $\hat{T}'$. The algorithm tries to prove by contradiction on $\mathcal{A}'$ the impossibility of each sub-branch. If the branch is proven impossible, we know (at confidence level at least $c$) that it cannot lead to an optimal solution, and it is pruned. Otherwise, it is added to the stack to be iterated on later. If all sub-branches are discarded (not represented in Fig. \ref{fig:flow}), we flag $x_{k+1}$ as an outlier and add $\mathcal{X}_{(x_1,\ldots,x_k,-1)}$ to the stack. This procedure leads to a small number of leaves that represent possible optimal solutions for $\boldsymbol{x}$.
 \paragraph{Generation of sub-branches}
 Using a prediction interval for $y_{k+1}$ based on $\hat{T}$, one can compute back possible values for $x_{k+1}$ at confidence level $c$. This problem is known as linear calibration (detailed in section 6.1.5 of \cite{Seber_Lee_2003}).
 Set $\bar{x}$ the arithmetic mean of the first $k$ items of $\boldsymbol{x}$ and $SS_{\boldsymbol{x}} = \sum_{i=1}^k\left(x_i-\bar{x}\right)^2$. By OLS:
$\hat{T} \sim \mathcal{N}\left(T,\frac{\sigma^2}{SS_{\boldsymbol{x}}}\right)$.
We get the prediction interval \cite{Predictive} at confidence level $c$ for an observation $y_0$ of predictor $x_0$: 
\begin{align*}
y_0 &\in [\widehat{y}_0-\Delta_\sigma, \widehat{y}_0+\Delta_\sigma] \\
\text{with } \Delta_\sigma &= z_{1-\frac{c}{2}} \sigma \sqrt{1+\frac{1}{n}+\frac{\left(x_0-\bar{x}\right)^2}{SS_{\boldsymbol{x}}}}
\end{align*}
If the algorithm is not executed online, 
%the total number of points in $\boldsymbol{x}$ is known and% 
a Bonferroni correction can be applied to ensure the significance level $c$ across all iterations of the algorithm. Replace $\sigma$ (unknown) by its (known) upper-bound $\overline{\sigma}$. 
%By using our upper-bound $\overline{\sigma}$ instead of the sample variance and $t$-statistic, we leverage our extra knowledge to build a tighter confidence interval, crucial to the efficiency of our algorithm.
From this estimate, we compute the so-called discrimination interval delimited by $\barbelow{x}_{k+1}$ and $\bar{x}_{k+1}$, respectively lower and upper bound for $x_{k+1}$ at significance level $c$ if they exist, by solving for $x_0$ the inequality: 
\begin{align*}
\left(y_0-\bar{y}-\hat{T}(x_0-\bar{x})\right)^2 \leq \chi^2_{1,1-\frac{c}{2}} \overline{\sigma}^2\left(1+\frac{1}{n}+\frac{(x_0-\bar{x})^2}{SS_{\boldsymbol{x}}}\right)
\end{align*}
We then take the intersection of this interval with $\{x_k+1,\ldots,x_k+N\}$ to keep only valid values for $x_{k+1}$.

 \paragraph{Plausibility of a branch}\label{tests}
The efficiency of our algorithm relies on how quickly it can identify branches that are unlikely to lead to the solution. To this aim, a variety of tests can be used assuming $\mathcal{A}'$ and performed at step 5. in Fig. \ref{fig:flow}:
 \begin{itemize}
     \item Local errors: For all $1 \leq i \leq k+1$, check that $y_i$ lies in the updated prediction interval for $\hat{y}_i$ at confidence $c$ (Bonferroni-corrected for the number of points).
     % Following \cite{Predictive}, we have that for any $1 \leq i \leq k+1$, 
     % $\hat{y}_{i} \sim \mathcal{N}(\overline{y_k},\sigma(1+1/n)
     % )$, with $\overline{y_k}$ the empirical mean of the $k$ first components of $y$.
     % We can thus once again take advantage of our upper-bound $\bar{\sigma}$ to compute a prediction interval at significance level at least $c$ (potentially Bonferroni-corrected) for each $\hat{y}_i$ and check that $y_i$ lies inside. If not, given our assumption for $\boldsymbol{x}$, the probability that this situation happens is less than $c$, and we can thus discard the branch.
     \item Global error: We have $r = \sum_{i=1}^n (\hat{y}_i-y_i)^2\sim \sigma^2\chi_{n-2}^2$. Using $\overline{\sigma}$, compute a confidence interval at significance $c$ and check that $r$ lies inside. Similarly, the mean absolute error should follow the mean of half-normal distributions.
     \item Normality tests: Run a Shapiro-Wilk or a K-S test at significance $c$ to check whether the distribution of residuals is indeed a multivariate Gaussian. This test is the only one having a super-linear complexity of $O(nlog(n))$ and can thus be skipped if asymptotic performance is an issue.
     \item Outliers: Assuming an upper bound $\overline{p}$ for the frequency of outliers, prune the branch if too many points are flagged as such. If the outliers are assumed uniformly distributed, check $1-G(n\_outliers) \leq c$ with $G$ the cumulative distribution function of a binomial law $\mathcal{B}(\overline{p},k+1)$.
     \item Outliers in a row: Check if too many consecutive points are flagged as outliers. A threshold can be chosen, or assuming uniformity one tests whether $\overline{p}^{o}\geq c$ with $o$ the current number of outliers in a row.
     \item Bounds on problem's variables: One can specify a range of acceptable values for $T$. If its intersection with the confidence interval from $\hat{T}'$ is empty, discard the branch. %Likewise, an arbitrary upper bound $N$ for the distance between two pulses can be used.
 \end{itemize}
 \paragraph{Indiscernible outliers}
  It is worth mentioning that if an outlier $(y_o,x_o)$ is not out of distribution, it will be integrated as a true pulse. For example, if pulse $10$ is missing but $y_o$ is indiscernible from $10T+b$ at confidence $c$, the algorithm will consider $x_o = 10$. While this has little impact on the result, it theoretically prevents our estimator from reaching the CRLB.
 \paragraph{Stopping at first leaf and optimality of the solution}\label{optimality}
 The tree traversal leads to a small number of leaves representing possible optimal solutions for $\boldsymbol{x}$. In practice, one can stop the algorithm as soon as a first leaf $\boldsymbol{x}^{\ast}$ is reached. Indeed, recall that we assumed $\sigma$ to be small with respect to $T$ so that each pulse cannot be interpreted as belonging to another period because of the noise. For $\boldsymbol{x}^{\ast}$ to pass the final local error tests we thus need $\forall i \in I, (y_i-b^\ast)/T^\ast - \lfloor (y_i-b^\ast)/T^\ast \rfloor \approx 0$. This implies that $\boldsymbol{x}^\ast$ is an approximate multiple of $\boldsymbol{x}$. By the integrality of $\boldsymbol{x}$, without indiscernible outliers, we thus have $\exists \lambda \in \mathbb{N}^*$ such that  $\boldsymbol{x}^\ast = \lambda\boldsymbol{x}$. Since we use a depth-first search, the leaf corresponding to $\lambda = 1$ will be reached first. To account for outliers, we check if $\exists \lambda \neq 1$ such that $\frac{1}{n}(\sum_{i=1}^n x_i^\ast \  \texttt{mod} \ \lambda) < \overline{p}$ and output $\lambda T^\ast$ if so.
 \paragraph{linear reconstruction}
%After a certain number of points $m$, the length of the discrimination interval becomes less than $1/2$, meaning that one can undoubtedly compute $x_{k+1}$ from $y_{k+1}$ and $\hat{T}$. In our tests, this usually happened between five to a few dozen points. This observation allows us to compute the amortized complexity of $O(\sum_{k}^n k) = O(n^2)$ for the algorithm (if only linear tests are used).
In the solution branch, $\hat{T}$ eventually becomes precise enough to undoubtedly compute (discrimination interval less than $1/2$) any $x_{k+1}$ from $x_k$ and $y_{k+1}-y_k$. In our tests, this usually happened between five to a few dozen points. At that time, a solution is found and one can reconstruct iteratively any remaining $x_i$, only checking if it is an outlier by the local error test. This way of ending the algorithm, called $\texttt{linear reconstruction}$ in Algorithm \ref{alg:iteration}, enables an amortized complexity of $O(n)$.
 
 \paragraph{Outliers in initial points}\label{assumptions}
 We assumed the first two points not to be outliers to have a starting point for the tree, and to be able to have a first estimation of $\hat{T}$. If these two assumptions are not met, the algorithm won't provide any solution. In that case, we rerun the algorithm without the first points. We chose to cap the number of such re-executions to 2. The code for the initialization is detailed in Algorithm \ref{alg:init}.
 
%as an iteration on branch $\mathcal{X}_(x_1,\ldots,x_k)$ does not make use of $(x_{k+1},\ldots,x_n)$ and the wrong branches are usually identified with only a couple of points
%%%% plausible
% \begin{algorithm}
% \caption{Plausibility checker}\label{alg:init}
% \begin{algorithmic}[1]
% \Procedure{is\_plausible}{$y,(x_0,\ldots,x_{i+1}),\theta$}
% \State $\hat{T},\hat{b} \gets \texttt{linear\_regression}(y,(x_0,\ldots,x_{i+1}))$
% \State $[\barbelow{T},\bar{T}] \gets \texttt{confidence\_interval}(c,\sigma,\hat{T},(x_0,\ldots,x_i))$
% \State $\sigma \gets \min(\bar{\sigma},\texttt{sample\_variance}(x_0,\ldots,x_{i+1})$
% \Ifm
%     \Sate $\exists 0 \leq k \leq i+1, |y_i - \hat{T}x_k - b|>\texttt{inverse\_CDF\_gaussian}(1-c/2,\sigma) \textbf{ or}$
%     \State $\texttt{MAE}(y - (\hat{T}\hat{x} + \hat{b}))>\texttt{inverse\_CDF\_half\_gaussian}(1-c/2,\sigma)  \textbf{ or}$
%     \State $\texttt{not }\texttt{Kolmogorov\_Smirnov}(c,y,(x_0,\ldots,x_{i+1}))   \textbf{ or}$
%     \State $1 - \texttt{CDF\_binomial}(\vars{n\_outliers},proba\_outlier,i+1) \geq c    \textbf{ or}$
%     \State $\vars{outliers\_in\_row} \geq \texttt{ceil}(\log(c)/\log(proba\_outlier))    \textbf{ or}$
%     \State $\vars{\bar{T}} > \vars{upper\_bound\_T}    \textbf{ or}$
%     \State $\vars{\barbelow{T}} > \vars{lower\_bound\_T}$
% \Then
%     \State $\Return \texttt{ False}$
% \EndThen
% \State $\Return \texttt{ True}$
% \EndProcedure
% \end{algorithmic}
% \end{algorithm}

\section{Tests}\label{section:test}
%%%%%%%%%%%%%%%%%%%%%%%%%%%%%%%% Init
\begin{algorithm}
\caption{Initialization}\label{alg:init}
\begin{algorithmic}[1]
\State $c \gets$ chosen confidence level
\State $\overline{\sigma} \gets $ best known upper bound for variance of $\epsilon$
\State $\barbelow{T},\overline{T}, N \gets$ arbitrary/physical bounds for $T$ and $x_{i+1}-x_i$
\State $\overline{p} \gets$ upper bound probability outlier
\State $\theta \gets (c,\overline{\sigma},N,\overline{p},\barbelow{T},\overline{T})$ %// all problem-dependant parameters
\Procedure{initialization}{$\boldsymbol{y},\theta$}
\State$\boldsymbol{y}\gets\boldsymbol{y}- y_1$
\State $\hat{T},\hat{b},x_1,x_2,\bar{x}_2 \gets 0,0,0,1,\min(N,\texttt{floor}(y_2/\barbelow{T}))$
%\State $\bar{x}_2 = \min(N,\texttt{floor}(y_2/\barbelow{T}))$
\While{$x_{2} \leq \bar{x}_{2} \texttt{ and not } \hat{T}$}
    \If{$\texttt{is\_plausible}(\boldsymbol{y},(x_1,x_2),\theta)$} 
        \State $T^\ast,b^\ast,x^\ast \gets \texttt{NextPoint}(\boldsymbol{y},(x_1,x_2),\theta)$
    \EndIf
    \State $x_2 \gets x_2 +1$
\EndWhile
\If{$\hat{T}$}
\State $\lambda = \texttt{biggest\_approx\_divisor}(x^\ast,\overline{p})$
\State \Return $\lambda T^\ast$
\EndIf
\State \Return $\texttt{initialization}((y_5,\ldots,y_n),\theta)$
\EndProcedure
\end{algorithmic}
\end{algorithm}

\begin{algorithm}
\caption{Branch iteration}\label{alg:iteration}
\begin{algorithmic}[1]
\Procedure{NextPoint}{$\boldsymbol{y},(x_1,\ldots,x_k),\theta$}
\State $outlier \gets \texttt{True}$
\State $\hat{T},\hat{b} \gets \texttt{linear\_regression}(\boldsymbol{y},(x_1,\ldots,x_k))$
\If{$\texttt{len}(\boldsymbol{x})=\texttt{len}(\boldsymbol{y})$}
\State \Return $\hat{T},\hat{b},\boldsymbol{x}$
\EndIf
\If{$\texttt{precise\_enough}(\boldsymbol{x},\boldsymbol{y},\hat{T})$}
\State \Return $\texttt{linear\_reconstruction}(\boldsymbol{x},\boldsymbol{y})$
\EndIf
\State $\barbelow{x}_{k+1},\bar{x}_{k+1} \gets \texttt{linear\_calibration}(\boldsymbol{y},\boldsymbol{x},\hat{T})$
\State $\barbelow{x}_{k+1},\bar{x}_{k+1} \gets \max(\barbelow{x}_{k+1},x_k+1),\min(\bar{x}_{k+1},x_k+N)$
%\State $\barbelow{x}_{i+1},\bar{x}_{i+1} \gets \texttt{linear\_calibration}(\boldsymbol{y},\boldsymbol{x},\hat{T})$
\For{$\barbelow{x}_{k+1} \leq x_{k+1} \leq \bar{x}_{k+1}$}
    \If{$\texttt{is\_plausible}(\boldsymbol{y},(x_1,\ldots,x_{k+1}),\theta)$}
        \State $outlier \gets \texttt{False}$ 
        \State \Return $\texttt{NextPoint}(\boldsymbol{y},(x_1,\ldots,x_k,x_{k+1}),\theta)$
    \EndIf
\EndFor
\If{$outlier$}
\State $\texttt{remove\_value\_from\_list}(y_k,\boldsymbol{y})$
\State \Return $\texttt{NextPoint}(\boldsymbol{y},(x_1,\ldots,x_k),\theta)$
\EndIf
\EndProcedure
\end{algorithmic}
\end{algorithm}
%%%%%%%%%%%%%%%%%%%%%%%%%%%%%%%%%%%%%%%

\begin{figure}
    \centering
    \includegraphics[width=\columnwidth]
        {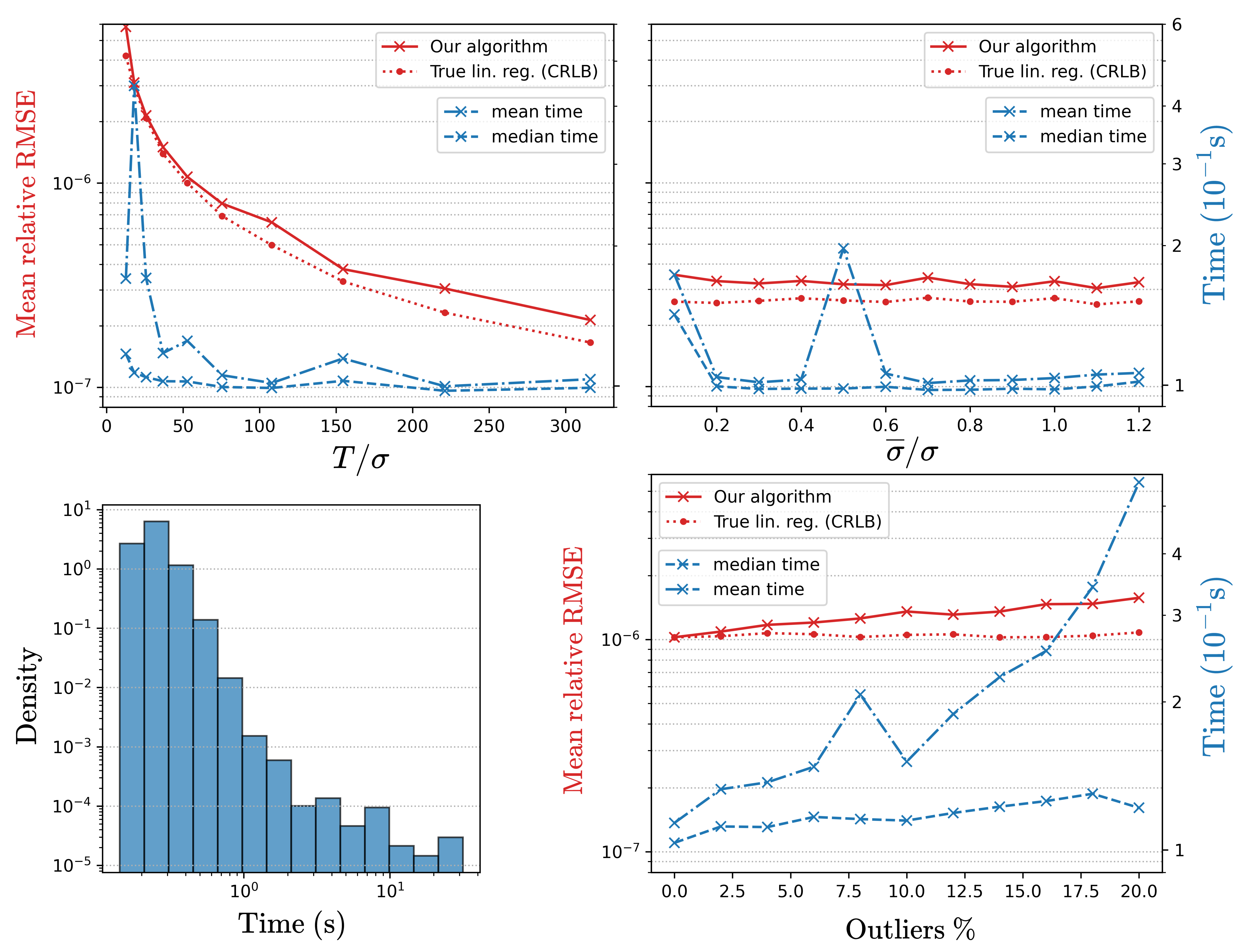}
    \caption{Duration and relative RMSE against different parameters. Each point is the mean result of 1,000 tests. The algorithm is executed 
    with $\overline{p}=0.1$ for \texttt{outliers\%} less than 0.1 and $\overline{p}=0.2$ otherwise.
    }
    \label{fig:tests}
\end{figure}
We first tested our algorithm on pulse indices generated by a Bernoulli process $\mathcal{B}(\gamma,1000)$ with sparsity $\gamma\sim \mathcal{U}(0.2,1)$. The histogram of Fig. \ref{fig:tests} displays the running times on a modern laptop as of 2023 for 10,000 signals of parameters:
\begin{itemize}
    \item $\sigma = 1$, $T \sim \mathcal{U}[20,100]$, $b \sim \mathcal{U}[0,T]$
    \item $\mathcal{U}\{0,\ldots,0.05n\}$ outliers uniformly sampled in $[0,y_n]$.
    \item $c = 10^{-6}$, $\overline{\sigma} = 1$, $N=50$.
\end{itemize} 
Most of the time, the algorithm is fast, with a median running time of $0.23$ seconds. However, in rare cases (notice the log-log scale), the execution degenerates to tens of seconds. This situation arises when some outliers are present among the very first points and make the tree grow deeper before being pruned, or force the algorithm to restart.   

It is worth mentioning that each branch is independent from the other ones, making this algorithm particularly suitable for parallelization. Moreover, it may be executed online. If necessary, its running time can thus be significantly reduced.

\begin{figure}
    \centering
    \includegraphics[width=\columnwidth]{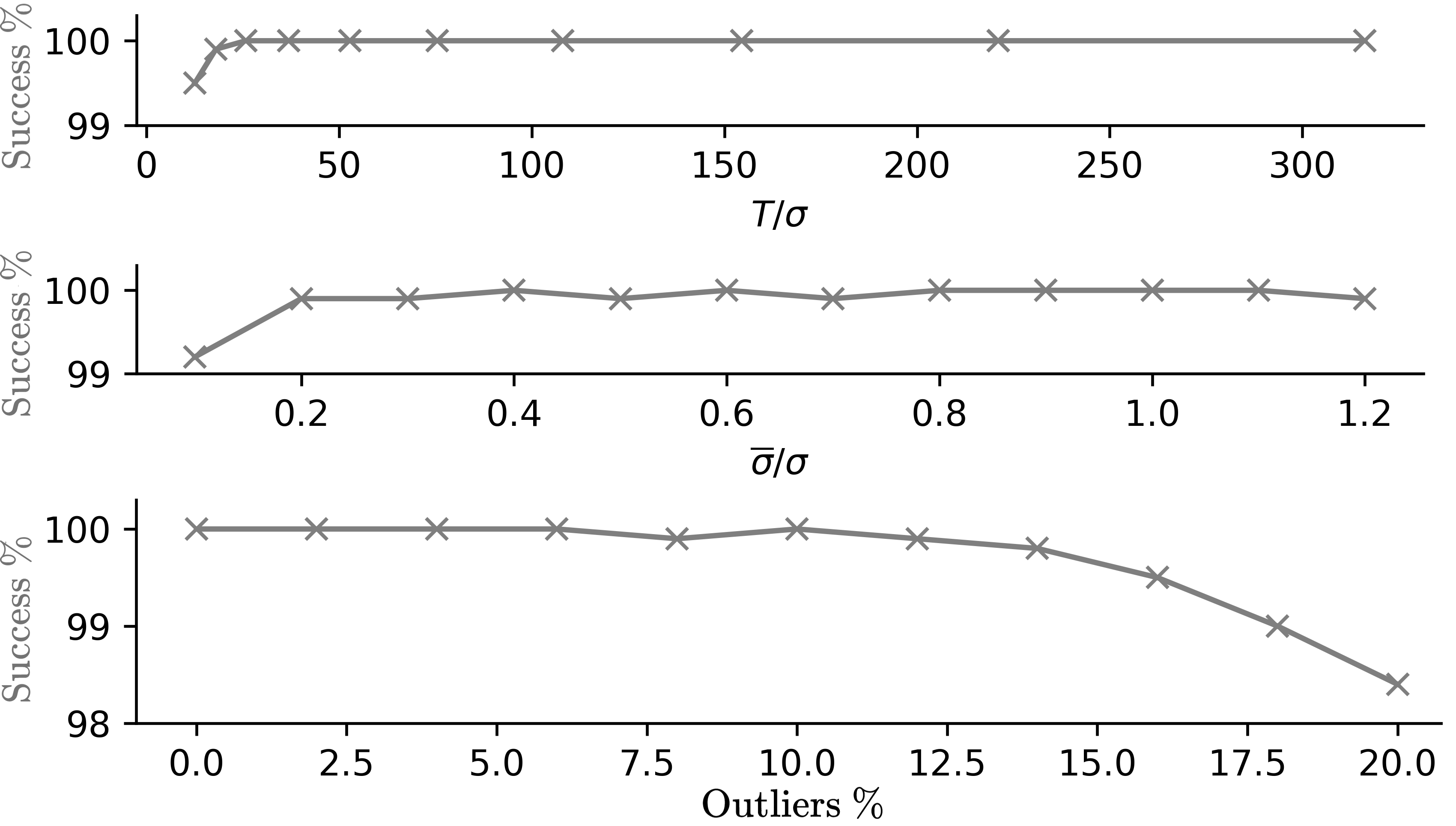}
    \caption{Success rate against parameters of Fig \ref{fig:tests}. An execution ``failed" if it did not output a solution, if its runtime exceeded one minute, or if the true solution is an integer multiple of the output one.}
    \label{fig:failures}
\end{figure}

The other panels of Fig. \ref{fig:tests} display the precision and running time of our algorithm against various parameters. Fig \ref{fig:failures} displays its success rate for the same parameters. The bottom-right panel of Fig. \ref{fig:tests} assesses the algorithm's efficiency against the quality of our upper bound $\overline{\sigma}$. The ratio $\sigma/T$ is fixed to $1/200$ to allow for large values of $\overline{\sigma}$ while still meeting our assumption $\overline{\sigma}/T < 1/10$ to preserve the pulses' order in time. Our simulations show that it is not necessary to provide an excellent bound on $\sigma$ for the algorithm to work efficiently, as long as the period-to-noise ratio $\sigma/T$ is low enough. The algorithm proves to be resilient in that regard, achieving near CRLB for ratios as high as $1/20$, close to the performance threshold of $\approx 1/10$ of the best estimators \cite{linear}. The upper-right panel demonstrates the state-of-the-art robustness of the algorithm to outliers, yielding precise estimations with up to 20\% contamination, even for high sparsity $\gamma$. In all tests, our algorithm achieves a precision of the order of the CRLB, and reaches it in the absence of outliers.

\section{Conclusion}\label{section:conclusion}
We have considered the problem of estimating the period $T$ of
a signal from sparse and noisy TOAs, contaminated by outliers, assuming that a decent upper bound on $\sigma$, the variance of the noise, is known. We provide an online, parallelizable, optimization-based algorithm to recover the pulses' indices and find the period by OLS, reaching the CRLB in the absence of outliers. In contrast to existing methods \cite{improved_quasi_ml}, our algorithm handles outliers by design with minimal effect on precision even with contamination up to 20\%, for various pulses' sparsity \cite{periodic_pulse_signal_analysis} and noise level. Our algorithm recovered near-CRLB solutions for ratios $\sigma/T$ as high as $1/20$, including up to 5\% of outliers, making it remarkably resilient with respect to the state-of-the-art \cite{linear}. Its complexity is amortized to $O(n)$, on par with the fastest estimators \cite{linear}. 
  
%A natural extension of this work would be to consider a multiple linear regression model, representing nested periods in the signal.

%This algorithm finds potential applications in situations where the signal is particularly noisy or when a very fine precision is needed.
%\paragraph{Generalization to multiple linear regression}The 2-D and n-D cases
%A natural and mathematically appealing extension of this problem is the generalisation to higher dimensions. One can for example consider signals that have repreating pattern, further subdivided in shorter periods. 

\bibliography{references}

\end{document}